\documentclass{article}
\usepackage[english]{babel}

\usepackage[letterpaper,top=2cm,bottom=2cm,left=3cm,right=3cm,marginparwidth=1.75cm]{geometry}

\usepackage{amsmath,amsthm,amssymb,mathtools}
\usepackage{bbm} 
\usepackage{graphicx}
\usepackage{subcaption}
\usepackage[colorlinks=true, allcolors=blue]{hyperref}
\usepackage[dvipsnames]{xcolor}

\usepackage{authblk} 

\def\R{{\mathbb R}}

\theoremstyle{plain}
\newtheorem{theorem}{Theorem}[section]
\newtheorem{proposition}[theorem]{Proposition}

\theoremstyle{definition}

\newcommand{\bv}{\overline{v}}

\newcommand{\bu}{\overline{u}}

\newcommand{\bvphi}{\overline{\varphi}}

\newcommand{\mmu}{\boldsymbol{\mu}}

\newcommand{\nn}{\mathbf{n}}

\newcommand{\xx}{\mathbf{x}}

\newcommand{\vK}{v_{K}}

\newcommand{\vL}{v_{L}}

\def\escalar#1#2{\left(#1,#2\right)}
\def\escalarL#1#2{\escalar{#1}{#2}}

\def\salto#1{\left[\!\left[#1\right]\!\right]}
\def\media#1{\left\{\!\!\left\{#1\right\}\!\!\right\}}

\def\T{\mathcal{T}}
\def\E{\mathcal{E}}

\def\O{\mathcal{O}}
\def\SO{\sigma}
\def\X{\chi}

\def\Pd{\mathbb{P}^{\text{disc}}}
\def\Pc{\mathbb{P}^{\text{cont}}}

\def\aupw#1#2#3{a_h^{\text{upw}}(#1;#2,#3)}

\title{Mathematical modeling of Neuroblast Migration towards the Olfactory Bulb}

\author[1]{Daniel Acosta-Soba}
\author[3,6]{Carmen Castro}
\author[4,6]{Noelia Geribaldi-Doldán}
\author[2,7]{Francisco Guillén-González}
\author[5,6]{Pedro Nunez-Abades}
\author[1]{Noelia Ortega-Román}
\author[6]{Patricia Pérez-García}
\author[1]{J. Rafael Rodríguez-Galván}

\affil[1]{Departamento de Matemáticas, Facultad de Ciencias, Universidad de Cádiz}
\affil[2]{Departamento de Ecuaciones Diferenciales y Análisis Numérico, Facultad de Matemáticas, Universidad de Sevilla}
\affil[3]{Área de Fisiología, Facultad de Medicina, Universidad de Cádiz}
\affil[4]{Departamento de Anatomía y Embriología Humanas. Facultad de Medicina. Universidad de Cádiz}
\affil[5]{Departamento de Fisiología. Facultad de Farmacia. Universidad de Sevilla}
\affil[6]{Instituto de Investigación e Innovación Biomédica de Cádiz (INIBICA)}
\affil[7]{Instituto de Matemáticas de la Universidad de Sevilla (IMUS)}

\newcommand{\Olf}{\mathcal{O}}
\newcommand{\SVZdomain}{\text{\it SVZ}}
\newcommand{\NZdomain}{\text{\it NZ}}
\newcommand{\CCdomain}{{C_C}} 
\newcommand{\CCexp}{\mu_{O}}

\begin{document}
\maketitle

\begin{abstract}
This article is devoted to the mathematical modeling of migration of neuroblasts, precursor cells of neurons, along the pathway they usually follow before maturing. This pathway is determined mainly by attraction forces, to the olfactory bulb, and the heterogeneous mobility of neuroblasts in different regions of the brain. In numerical simulations, the application of novel discontinuous Galerkin methods allows to maintain the properties of the continuous model such as the maximum principle. We present some successful computer tests including parameter adjustment to fit real data from rodent brains.
\end{abstract}

\section{Introduction}

Adult neurogenesis has garnered significant attention in recent years with many authors trying to shed light into this relevant process.
New neurons are generated in the adult rodent brain in specialized regions in which  neural stem cells (NSC) are activated to produce neurons in a hierarchical process named neurogenesis. One of the main neurogenic regions in the adult rodent brain  is the subventricular zone (SVZ)~\cite{obernier_2019}. Activation of SVZ NSC induces their cell cycle entrance and posterior division to produce transit amplifying actively dividing progenitors that will give rise to neuronal progenitors or neuroblasts ~\cite{Codega_2014}. Newly generated neuroblasts migrate from the SVZ toward the olfactory bulb (OB) through a path called rostral migratory stream (RMS) while still dividing ~\cite{Ponti_2013}, contributing to the continuous neuronal replacement in the OB. Neurons produced in this homeostatic mechanisms integrate into existing circuits and participate in olfaction ~\cite{Lazarini_2011}. 

Tens of thousands of neuroblasts migrate daily through the RMS travelling long distances toward the OB, where they integrate as inhibitory interneurons~\cite{Lois_and_Alvarez-buylla_1994, altman,Morshead,carleton_lledo_2003}. These cells apparently display an exploratory pattern while moving long distances ~\cite{Nam}. On this intrincated way along the adult RMS, neuroblasts border the corpus callosum (CC), a bundle of nerve fibers located between the left and right cerebral hemispheres, where the movement of neuroblasts is hindered due to the presence of a glial sheath delineating the RMS and the low density of blood vessels in the CC, which they use as scaffold ~\cite{Monyer_2012_migration_vessels}.

The importance of the homeostatic migration of neuroblasts toward the OB in the adult brain is highlighted by the fact that in models of brain damage an altered migration pattern is found. Several models of brain damage show an altered migration of SVZ neuroblasts that are different depending on the type of damage. Thus, in a murine model of Alzheimer's disease  the proportion of migrating cells in the RMS is lower than in healthy mice ~\cite{Esteve_RMS_AD} whereas in other models, such as murine models of Huntington's disease migration toward the OB is reduced as well and  accompanied by an altered pattern of migration that direct neuroblasts toward the striatum avoiding the OB pathway ~\cite{Kohl_2010_migration_OB_huntington, Kandasamy_migration_huntington_rat}. Alterations of neuroblast migration are also observed in cortical injuries generated by ischemic lesions, in which SVZ neurogenesis is stimulated in response to the injury and in some cases chains of neuroblast can be seen migrating toward the injured region~\cite{magavi_migration_cortex}. 
All these reports suggest that stimulating  neurogenesis in the SVZ and conducting neuroblasts toward the injured region may be of use at designing strategies to regenerate damaged brain regions. 

Notwithstanding, although neurogenesis and neuroblast migration in the adult brain of mammals has been studied in depth over the past two decades from an anatomical and physiological point of view, and despite the considerable amount of data available over the years, up to our knowledge no mathematical or computational model has been published to date describing the neuroblast migration toward the OB in realistic brain domains. 
And this type of models could be of great value as a first approximation to describe more complex phenomena, such as the movement of neuroblasts towards brain lesions.

Virtually all the deterministic neurogenesis mathematical models in the literature are based on compartmental systems of ordinary differential equations (ODEs), see e.g.~\cite{ziebell2014mathematical, ziebell2018revealing, kalamakis2019quiescence, dabelow2022distinguishing} or the recent review~\cite{danciu2023review}. 
They delineate, over time, the asymptotic population dynamics of various cell states (such as NSCs, neuroblasts, neurons, etc.) that play hierarchical roles in the neurogenesis process. These models also incorporate compartments representing certain properties of these cells, such as quiescence, proliferation, self-renewal, apoptosis, among others. 

These ODE models do not cover spatial movement of neouroblasts, hence more general partial differential equations (PDE) models would be required for that. To our knowledge, the only interesting contribution in this direction are the paper by Ashbourn et al.~\cite{Ashbourn2012mathematical} and its subsequent expansion in~\cite{VanSchepdael2013mechanisms}. In these references, a few of the key ideas in which we are independently basing this paper can be found. Notably, their model postulates a singular chemoattractant governing cell migration, emanating from both the olfactory bulb (OB) and also from the cells involved in neurogenesis. However, it is worth noting that PDEs are numerically solved within a considerably simpler and less realistic spatial domain compared to the framework proposed in our study. The algorithm used in these references to solve numerically the neuroblast PDEs is the finite-volume one presented in~\cite{gerisch2006robust}. In particular, modeling and design of a complex chemoattractant function driving migration of neuroblasts along a spatially complex RMS path is not covered in these works. 
 
In this paper we present a PDE mathematical model of neuroblast migration to the olfactory bulb, based on the following hypothesis: the movement of neuroblasts is due to a transport phenomenon exerted by certain attraction or \textit{chemotaxis} velocity towards the OB.
More specifically, we suppose that there exists a function to be identified whose gradient drive the transport of neuroblasts. 
One of our main tasks in this work will be the correct modeling and numerical computing of this olfactory bulb function, which is defined along the brain and whose isolines reproduce the RMS path around the CC.

Once identified the transport velocity, we focus on a second task: computing the migration of neuroblasts along the potentially realistic brain domain as solution of a convection PDE. It is well known that computing numerical solutions to this kind of hyperbolic PDE equations is not easy and numerical schemes must be carefully developed in order to avoid nonphysical oscillations and maintain desirable properties of the continuous model like the pointwise bounds~\cite{ern_guermond_2010,quarteroni2010numerical}. Fortunately we can successfully apply the upwind Discontinuous Galerkin (DG) techniques published in \cite{acosta-soba_upwind_2022} to obtain quite appropriate results.
The PDE is supplemented with some extra reaction terms described below, modeling the born of neuroblasts and their disappearance due to evolution in mature neurons. 

Altough technically not of classical chemotaxis type, our model is inspired by chemotatic biological processes through which a population of organisms (or cells) migrate in response to a chemical stimulus. More specifically the study of cell motion affected by a chemical gradient, which results in net propagation up a chemoattractant gradient or down a chemorepellent gradient.
Since the classical mathematical model, introduced in the 70's by E.F.~Keller and L.A.~Segel~\cite{keller_segel_1970,keller_segel_1971} the topic has aroused considerable interest in the mathematical community and many related variants have been proposed (see e.g.~\cite{bellomo2015toward} for a review).

\begin{figure}
  \centering
  \includegraphics[width=0.7\linewidth]{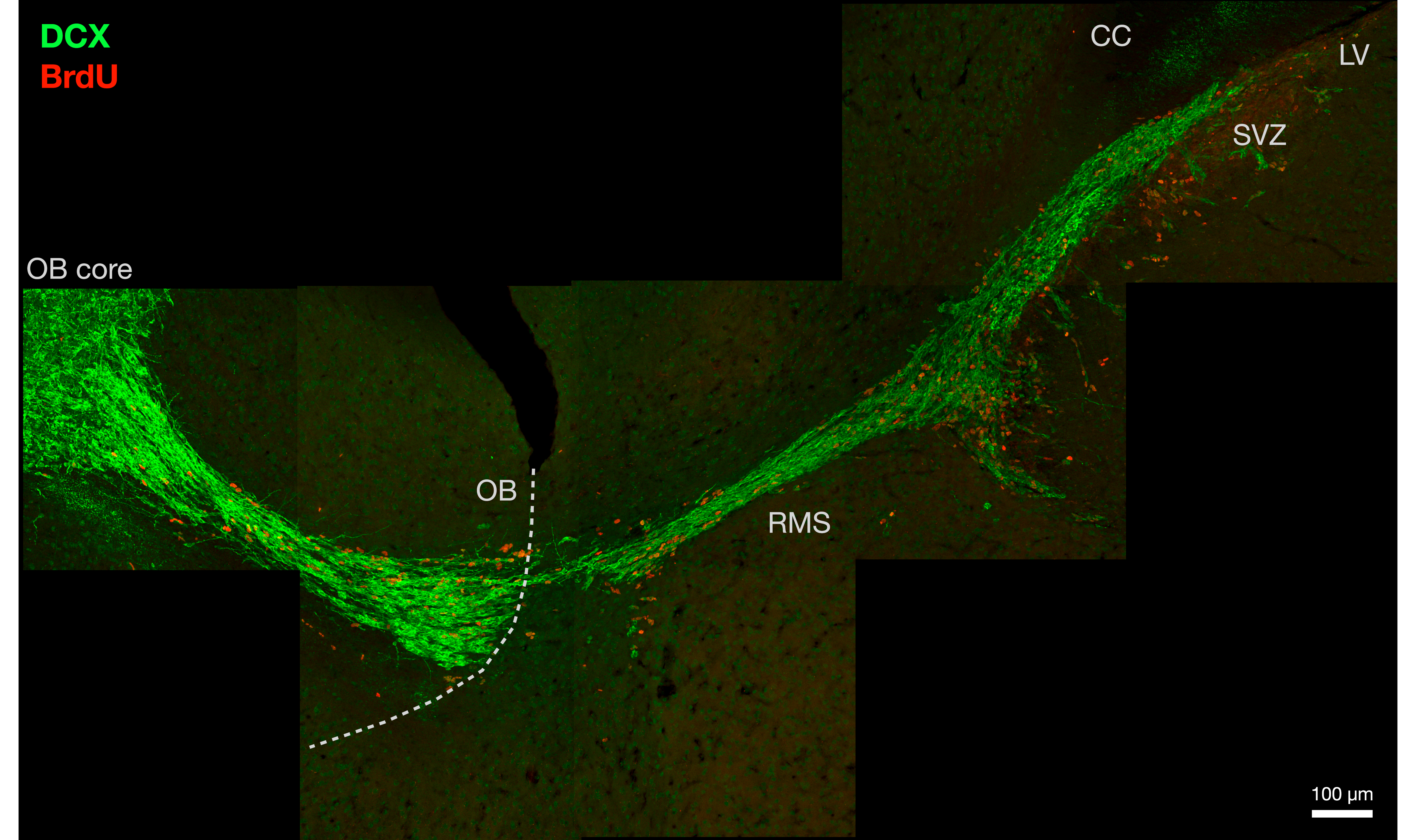}
  \caption{ Composition of confocal images showing neuroblasts migrating from the SVZ, to the OB through the RMS. To obtain this image, mice were given three intraperitoneal injections of the thymidine analogue BrdU in one day and were sacrificed 3 hours after the last BrdU injection. Brain were isolated and 30 µM sagittal sections were obtained with the use of a cryostat. Sections were then processed for immunohistochemistry to detect BrdU (red) and the neuroblast marker DCX (green). The scale bar represents 100 µm. Abbreviations: subventricular zone (SVZ); corpus callosum (CC); rostral migratory stream (RMS); olfactory bulo (OB)}
  \label{fig:neuroblast-real}
\end{figure}

Finally, it is worth remarking that our mathematical model for neuroblast migration not only has been implemented in the computer, but it also has been calibrated using real data that we have previously experimentally obtained. Specifically, we are using databases obtained from real images of the rodent brain, where the neuroblasts can be identified as the red spots in Figure~\ref{fig:neuroblast-real}. These images show specifically a 2D saggital section of the brain where neuroblasts are labelled by immunostaing with the neuroblast marker doublecortin (green) and with bromodeoxyuridine (BrdU) that was injected in the ventricle 3 hours  before sacrifice and was incorporated into SVZ dividing cells at the time of the injection. We intend to get some insights about the possibility of adjusting the parameters of our mathematical model so that the solution match not only qualitatively, but quantitatively the experimental counting in these real images of the brain.
In the end, we are exploring the possibility of obtaining a neuroblasts migration model in realistic brain domains which can be validated experimentally.

\section{The Model of Neuroblast Evolution}
\label{sec:neuroblast-model}

In this section, we go into details of the model for migration of neuroblasts towards the OB along the RMS, starting from the SVZ and rounding the CC. Our more relevant hypothesis are: (1) the shape of this RMS path can be modeled by means of an attraction function, which we are going to carefully design in the following sections, and (2) the migration of neuroblasts can be modeled by an adequate convection+reaction equation, where the convective velocity is related to the gradient of the previous attraction function.

\subsection{The Neuroblast PDE Model}

Going into details, let us fix a space-time domain $\Omega\times (0,T)$, where $\Omega\subset\R^2$, is an open set representing a rodent brain, with boundary $\partial\Omega$. We consider the evolution of the density of neuroblasts that have been marked with BrdU at the initial time $t=0$ by means of the following problem:
\begin{subequations}
\label{eq:neuroblast}
\begin{align}
\label{eq:neuroblast.eq}
\tau u_t + \X\, \nabla \cdot(u\nabla \O) + \alpha\, u - \gamma\, u \,\mathbbm{1}_{\NZdomain} 
= \beta\, \mathbbm{1}_{\SVZdomain}  
&\qquad\text{in}\enspace\Omega\times (0,T),
\\
\label{eq:neuroblast.bc}
 u=0
&\qquad\text{on}\enspace\partial\Omega\times (0,T),
\end{align}
\end{subequations}
where $\mathcal{O}=\mathcal{O}(\xx)\in \R$ is a potential function such that its gradient $\nabla \O$ models the attraction exerted by the olfactory bulb.  The precise definition of this function 
$\mathcal{O}$, reflecting the heterogeneity of the brain with respect to the neuroblast migration, is detailed in Section~\ref{sec:olfact-bulb-chem} below. The function $\mathcal{O}$ depends on a parameter $\SO>0$, related to some spatial distribution of the chemoattractant:  concentration around the center of the OB (small values of $\SO$) or spread along the domain (big values of $\SO$), in a Gaussian fashion related to the function~(\ref{source}). On the other hand, the parameter $\X>0$ is the amplitude of this attraction, $\alpha>0$ is the apoptosis or death rate along the domain and $\beta>0$ is related to the amount of neuroblasts which are generated in the SVZ, after a process of division of NSCs whose study will not be covered here.
Also a source term, with coefficient $\gamma>0$, has been added for the purpose of  modeling the contribution to our two-dimensional domain of neuroblasts coming from other regions of the real three-dimensional brain. In particular, we assume that this contribution is located  zone of the rodent brain where it becomes tighter, as can bee seen e.g. in~\cite{atlas}. This area will be called the narrowing zone (NZ). 

In equation~\eqref{eq:neuroblast}, we assume the following condition for the chemoattractant function $\mathcal{O}$: 
\begin{equation}
  \label{eq:nabla.O.dot.n}
\nabla  \mathcal{O}\cdot \mathbf{n} < 0 \enspace\text{on}\enspace \partial\Omega,
\end{equation}
where $\textbf{n}$ is the exterior unit normal vector $\partial\Omega$. 
Therefore the inflow boundary 
$$
\partial\Omega^- = \{ x\in\partial\Omega \ /\ \nabla\mathcal{O}\cdot\nn < 0 \}
$$
consist of the entire boundary $\partial\Omega$ and hence the Dirichlet boundary condition~(\ref{eq:neuroblast.bc}) is imposed.
Finally, the parameter $\tau\in\{0,1\}$ distinguishes two cases which are modeled by~(\ref{eq:neuroblast}): the stationary model (corresponding to $\tau=0$) and the evolution model ($\tau=1$).
We are going to focus mainly on the last one, $\tau=1$, equation~(\ref{eq:neuroblast}), in which case  an initial condition
\begin{equation}
  \label{eq:neuroblast.ic}
  u(0)=u^0 \enspace\text{in}\enspace \Omega,
\end{equation}
is added,
where $u^0$ is a function defining the initial density of neuroblasts along $\Omega$. 

In this work we are interested in calibrating the parameters of our evolution model~(\ref{eq:neuroblast}),
\begin{equation}
\label{eq:parameters}
    \Lambda = (\alpha, \beta, \gamma, \X, \SO)\in\mathbb{R}^5_+,
\end{equation}
testing if this model can can be fitted to experimental rodent brain data. Specifically, we would like to computationally reproduce the data set which we had previously obtained from microscopical images and tabulated at different time instants (see Figure~\ref{fig:neuroblast-real}, for a real microscopical image of the distribution of neuroblasts, and Tables~\ref{table:real.data.t0} and~\ref{table:data.evolution} for the real data). For that purpose, the  initial value, $u_0$, must be as close as possible to the data at the initial time $t_0$.
This initial value, $u_0$, will be defined as the solution to the stationary ($\tau=0$) neuroblast equation~(\ref{eq:neuroblast}), assuming the model is in a steady equilibrium state, for another adequate parameter set. 

Lastly, it should be noted that in the evolution case, $\tau=1$, it does not make sense to consider a source of marked neuroblasts in the SVZ because those neuroblasts that are born in the SVZ after the injection is given (initial time) have not received the BrdU marker and therefore are not object of our model. Therefore we are going to fix $\beta=0$ in this case.

\subsection{The Olfactory Bulb Chemoattractant Function}
\label{sec:olfact-bulb-chem}

The attraction due to the olfactory bulb is computed in terms of a function, $\O=\O(\xx) \ge 0$, whose gradients dictate the transport of neuroblasts according to~(\ref{eq:neuroblast}). In this section we will show how this function is obtained as the weak solution to a steady problem with reaction and anisotropic diffusion. Our aim is that this ``olfactory bulb function'' $\O$ reaches its maximum at some given point $(x_\Olf,y_\Olf)$, which stands for the center point of the olfactory bulb. Moreover $\O$ will decrease homogeneously along $\Omega$, except in a certain area, the corpus callosum (CC), where the gradient of $\O$ vanishes. This way we are going to model the fact that it is hard for neuroblasts to move within this region.

In order to obtain a function with the above-mentioned characteristics, we proceed as follows and let $\CCdomain\subset\Omega$ be an open set defining the CC and let $(x_\Olf, y_\Olf)\in\Omega\setminus\overline\CCdomain$. We define the following piecewise constant expression which will be used as an anisotropic diffusion coefficient in our olfactory bulb chemoattractant model:
\begin{equation}
\label{eq:OB-viscosity}
\CCexp =
\left\{
\begin{array}{ll}
\displaystyle
    \frac{1}{P_{C}} & \text{ in } \CCdomain,  \\ \\
     P_{C} & \text{ otherwise. }
\end{array}
\right.
\end{equation}
The parameter $P_{C}\in (0,1)$ is a small positive value, $P_{C} \ll 1$, which can be understood as the degree of ``permeability'' of the CC regarding  neuroblast migration across this region of the mouse brain. 
The idea is the following: to assign a huge diffusion number, $\CCexp=\frac{1}{P_{C}} \gg 0$, to $\CCdomain$ while an almost vanishing diffusion $\CCexp\simeq 0$  to the rest of the brain.

On the other hand, we consider the following source term, defined as a Gaussian function centered at the middle point of the olfactory bulb $(x_{\O}, y_{\O})$:
\begin{equation}
    \label{source}
    f_{\O}(x,y) = f_{\O}(\SO; x,y) 
    =  e^{-((x-x_O)^2 + (y-y_O)^2)/\SO^2}.
\end{equation}
The constant $\SO>0$ is one of the key parameters in our model, as it determines the spread of the attraction of the OB along the brain and therefore the migration speed of neuroblasts in each point of $\Omega$. It will be estimated from real data by the process described in the following sections.

With all the above, we consider the following problem: find $\O\in H^1(\Omega)$ -- the Sovolev space of $L^2(\Omega)$ functions with weak derivatives in $L^2(\Omega)$ -- as the solution to
\begin{equation} 
  \label{eq:OB}
  \left\{ \begin{aligned} \Olf-\nabla\cdot (\mu_{O}\,\nabla \Olf) &= f_{\O}  \quad  \text{ in } \Omega, \\
      \Olf &= f_{\O}  \quad \text{ on } \partial\Omega.
    \end{aligned}
  \right.
\end{equation}
It is well known that this elliptic equation has an unique solution which is at least in $C^0(\Omega)$, see e.g. 
\cite[Theorem 8.24]{gilbarg1977elliptic}.
It is expected that the huge weight of the coefficient $\CCexp$ in the corpus callosum CC will make the solution $\O$ to take a practically constant value in that area. Therefore $\nabla\O$ will be very small in this region, so that migration of neuroblasts, defined by equation~(\ref{eq:neuroblast}), is negligible inside the CC. Whereas in the remaining domain, the vanishing viscosity considered in~(\ref{eq:OB-viscosity}) makes the solution behave like the anisotropic Gaussian function $f_{\O}$ and thus the movement of neuroblasts acts according to the gradient of $f_{\O}$ towards the center of the OB.  

Notice that we are considering this attraction independent of  the time.
Our numerical solutions (Figure ~\ref{fig:steady-opt-ob}) confirm the expected behavior for $\O$ (see e.g. Figure~\ref{fig:steady-opt-ob}). This behavior is one of the keys of our model since, as we will show  $\nabla\O$ is able, by itself, of defining the complex RMS path along the neuroblast migration occurs (see e.g. Figure~\ref{fig:qualitative} below).

\section{Numerical Solution to the Neuroblast Model}
\label{sec:numerical-schemes}

\subsection{Notation}
Let $0=t_0<t_1<\cdots<t_M=T$ be an uniform partition of the time domain $[0,T]$ with $\Delta t=t_m-t_{m-1}$. Given any scalar function $v\colon \Omega\times[0,T]\longrightarrow\R$ and an approximation  $v^m\approx v(\cdot, t_m)$, we denote by   $\delta_t v^m=({v^m-v^{m-1}})/{\Delta t}$ the discrete backward time derivative operator.

For the space discretization, we consider a shape-regular triangular mesh $\T_h=\{K\}_{K\in \T_h}$ of size $h$ over  a bounded polygonal domain $\Omega\subset\R^d$. 
We define the approximation spaces of discontinuous, $\Pd_k(\T_h)$,  and continuous, $\Pc_k(\T_h)$, finite element functions over $\T_h$ whose restriction to $K\in\T_h$ are polynomials of degree $k\ge 0$.
Moreover, we note the set of edges or faces of $\T_h$ by $\E_h$, which can be split as $\E_h=\E_h^i\cup\E_h^b$ into the \textit{interior edges} $\E_h^i$ and the \textit{boundary edges} $\E_h^b$. Now, we fix the following orientation for the unit normal vector $\nn_e$ associated to any edge $e\in\E_h$ of the mesh $\T_h$:
 \begin{itemize}
 	\item Let $e\in\E_h^i$ be an interior edge with unit normal vector $\nn_e$. Then $e$ is shared by two elements of $\T_h$, which will denoted by $K_e$ and $L_e$, in such a way that $\nn_e$ is exterior to $K_e$ pointing to $L_e$.

If there is no ambiguity, in order to abbreviate the notation we will denote the previous elements $K_e$ and $L_e$ simply by $K$ and $L$, respectively, assuming that their naming is always with respect to the edge $e$ $\E_h^i$ and it may vary if we consider a different one.

 	\item If $e\in\E_h^b$, then $\nn_e$ points outwards of the domain $\Omega$.
 \end{itemize}



In addition, the \textit{average} $\media{\cdot}$ and the \textit{jump} $\salto{\cdot}$ of a scalar function $v$ on an edge $e\in\E_h$  are defined as follows:

\begin{equation*}
		\media{v}\coloneqq
		\begin{cases}
			\dfrac{\vK+\vL}{2}&\text{if } e\in\E_h^i\\
			\vK&\text{if }e\in\E_h^b
		\end{cases},
		\qquad
		\salto{v}\coloneqq
		\begin{cases}
			\vK-\vL&\text{if } e\in\E_h^i\\
			\vK&\text{if }e\in\E_h^b
		\end{cases}.
\end{equation*}

Finally, we set the following notation for the positive and negative parts of a scalar function $v$:
	$$
	v_\oplus\coloneqq\frac{|v|+v}{2}=\max\{v,0\},
	\quad
	v_\ominus\coloneqq\frac{|v|-v}{2}=-\min\{v,0\},
	\quad
	v=v_\oplus - v_\ominus.
	$$


\subsection{Computing the Olfactory Bulb Attraction}
\label{sec:computing.olf}
The linear reaction-diffusion equation~\eqref{eq:OB} in the domain defined by the mesh $\T_h$ will be solved by means of classical a $P_1$-continuous Finite Elements Method (FEM). For that, we have to deal with the discontinuity of the viscosity coefficient $\CCexp$, what is done as follows.
We project this coefficient into a discontinuous piecewise constant function  $\mu_{O,h} \in \Pd_0(\T_h)$, which is defined on each element $K$ of $\T_h$ as: $\mu_{O,h}|_K=1/{P_{C}}$ if the barycenter of $K$ is in the CC,  $\mu_{O,h}|_K={P_{C}}$ otherwise.
An example of this kind of piecewise constant approximation can be seen  e.g. in Figure~\ref{fig:brain-mesh}, where triangles related to the CC are colored in black. 

Once this done, given the parameter $\SO>0$, let  $f_\O$ be the Gaussian function introduced in~(\ref{source}). Let us consider the problem: find $\Olf\in\Pc_1(\T_h)$ such that $\Olf = f_{\O}$ on $\partial\Omega$  and 
\begin{equation}
  \label{esquema_DG_OlfBulb}
  b(\O, \bvphi)=\escalarL{f_{\O}}{\bvphi}, 
  \enspace
  \forall \bvphi\in V_h, 
\end{equation}
where $(\cdot,\cdot) = (\cdot,\cdot)_{L^2(\Omega)}$ stands for the scalar product in $L^2(\Omega)$, $V_h = \{\bvphi\in\Pc_1(\T_h), \ \bvphi|_{\partial\Omega}=0\}$ and 
$$
b(\O, \bvphi)\coloneqq 
\int_{\Omega}   \O \ \bvphi + \sum_{K\in T_h} \int_K \mu_{O,h} \nabla \O \cdot \nabla \bvphi   .
$$
%
%
The following result can be readly derived from the Lax-Milgram Theorem.
\begin{proposition}
  For any constant $\SO>0$, there is a unique solution to the problem \eqref{esquema_DG_OlfBulb}.
\end{proposition}

\begin{figure}
  \centering
  \includegraphics[width=.8\linewidth]{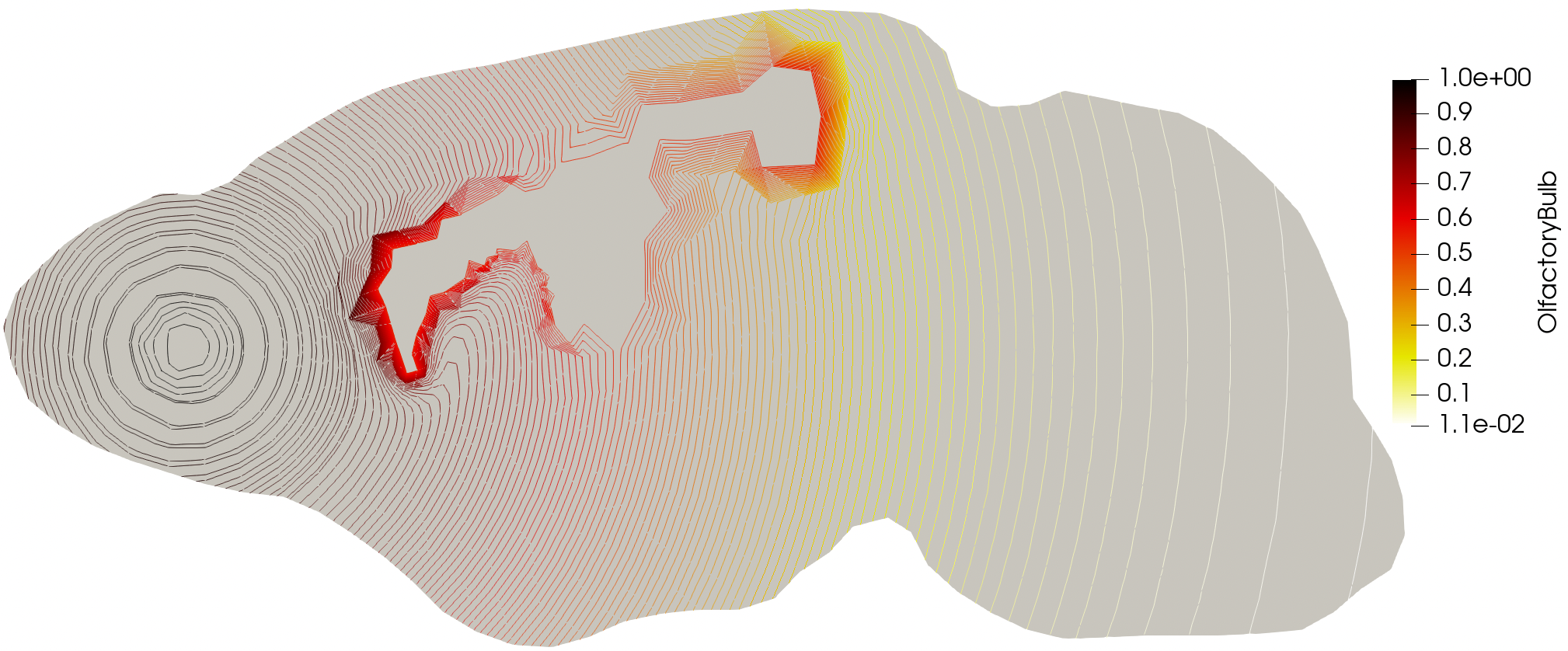}
  \caption{Isolines of the olfactory bulb attraction function, obtained as solution of~\eqref{esquema_DG_OlfBulb} with the parameters described in Section~\ref{sec:experimental-param}. The olfactory bulb stands out for the concentric isolines on the left. The corpus callosum is located in the center, in the region with no isolines.}
  \label{fig:steady-opt-ob}
\end{figure}

The resulting approximate solution $\O$ satisfies the properties indicated in Section~\ref{sec:olfact-bulb-chem}, in particular it is increasing towards the OB and $\nabla\O$ vanish in the CC, also $\nabla\mathcal{O}\cdot\nn\le0$ on $\partial\Omega$. For a concrete example, where $\sigma$ was computed from the process of adjustment of parameters described in Section~\ref{sec:experimental-param}, see Figure~\ref{fig:steady-opt-ob}.

\subsection{Computing the Initial Condition}
\label{sec:appr-init-cond}

At this point we can introduce the discrete scheme to approximate the solution to the neuroblasts migration model~\eqref{eq:neuroblast} for the stationary case $\tau=0$. This solution, $u_0$, will be taken as the starting point for the evolution model ($\tau=1$), see Section~\ref{sec:numer-solution}. 

The careful design of this discrete scheme is one of the key points in this work.
In fact, it is widely known (see e.g.~\cite{ern_guermond_2010, quarteroni2010numerical}) that a careless use of classical FEM methods for approximating convection-reaction equations like~\eqref{eq:neuroblast} can lead to spurious oscillations and awful results, related to the loss of the maximum principle property (MPP), which is verified at continuous level.
Here we will use a MPP-preserving upwind discontinuous Galerkin method introduced in~\cite{acosta-soba_upwind_2022}.

More in detail, for any set of parameters $\Lambda = (\X, \SO, \alpha, \beta, \gamma)\in\mathbb{R}^5_+$, we consider 
$\mmu\in \left[\Pc_1(\T_h)\right]^d$ as the $P^1$--continuous $L^2(\Omega)$ projection of $\nabla\O$, being $\O\in \Pc_1(\T_h)$ the solution to \eqref{esquema_DG_OlfBulb} for the parameter $\SO$ and define
the bilinear form
\begin{equation}
\label{eq:a_Lambda}
a_\Lambda(v,\bv) \coloneqq \X\, \aupw{\mmu}{v}{\bv} + \alpha\, \escalarL{v}{\bv} - \gamma\, \escalarL{v\mathbbm{1}_{\NZdomain}}{\bv},
\end{equation}
where the upwind bilinear form $\aupw{\cdot\;}{\cdot}{\cdot}$ is defined as follows (see~\cite{acosta-soba_upwind_2022}):
\begin{equation}\label{eq:aupw}
	\aupw{\mmu}{v}{\bv}
	= -\sum_{K\in\T_h}\int_K v\,\mmu\cdot\nabla\bv 
	\enspace + \sum_{e\in\E_h^i, e=K\cap L} \int_e\left( (\mmu\cdot\nn_e)_{\oplus}\vK-(\mmu\cdot\nn_e)_{\ominus}\vL\right)\salto{\bv},
\end{equation}
for any functions $v,\bv$ in the broken Sobolev space $H^1(\T_h)$ .
In the case of $v,\bv\in\Pd_0(\T_h)$, since $\mmu\cdot\nn\le0$ on $\partial\Omega$, the bilinear form~\eqref{eq:aupw} reduces to 
\begin{equation}
  \label{eq:aupwFinal}
  \aupw{\mmu}{v}{\bv} = 
  \X\!\!\!\! \sum_{e\in\E_h^i, e=K\cap L} \int_e\left( (\mmu\cdot\nn_e)_{\oplus}\vK-(\mmu\cdot\nn_e)_{\ominus}\vL\right)\salto{\bv},
\end{equation}
where the boundary condition~(\ref{eq:neuroblast.bc}) has been weakly imposed.

Let us fix an appropriate parameter set which will be denoted by
$$
\Lambda_0=(\X_0, \SO_0, \alpha_0, \beta_0, \gamma_0)\in\mathbb{R}^5_+,
$$ 
and which can be chosen for instance as in Section~\ref{sec:experimental-param}. 
The following linear form on $\Pd_0(\T_h)$ is considered:
\begin{equation}
\label{eq:L_Lambda}
L_{\Lambda_0}(\bv)\coloneqq\beta\,\escalarL{\mathbbm{1}_{\SVZdomain}}{\bv}.
\end{equation}
Then we approximate the steady ($\tau=0$) solution $u_0$ to~(\ref{eq:neuroblast}) by a piecewise constant function defined as the solution to the following problem: find $u^0\in \Pd_0(\T_h)$ such that
\begin{equation}
  \label{esquema_DG_neuroblasts_steady}
  a_{\Lambda_0}(u^0,\bu)=L_{\Lambda_0}(\bu), \quad \forall \bu\in\Pd_0(\T_h).
\end{equation}



\subsection{Numerical Solution of the Evolution Model}
\label{sec:numer-solution}

Now, let us consider:
\begin{enumerate}
\item A chemoattractant discrete function
  $\O\in\Pd_1(\T_h)$, computed as the solution to the olfactory bulb problem \eqref{esquema_DG_OlfBulb}, depending on a given parameter $\SO = \SO^0$ previously computed and which we will keep unchanged in what follows.
\item A discrete initial condition $u^0\in\Pd_0(\T_h)$, computed using the previously computed chemoattractant function $\O$ and for a given parameter set $\Lambda_0$, which might be calibrated as detailed in Section~\ref{sec:experimental-param}.
\item A given set of parameters $\Lambda = (\X, \SO, \alpha, \beta, \gamma)\in\mathbb{R}^5_+$, where the following values have been fixed: firstly $\sigma=\sigma^0$, because the shape of the OB is a biological characteristic which doest not difer between evolution and transient cases, and secondly $\beta=0$, as stated in Section~\ref{sec:neuroblast-model}.
\end{enumerate}
Then we approximate the solution of the evolution ($\tau=1$) system~(\ref{eq:neuroblast})--(\ref{eq:neuroblast.ic}) using a first order implicit-explicit time scheme and an upwind DG space discretization as follows.
Let us start from $u^0\in\Pd_0(\T_h)$. For each $m\in \mathbb{N}$, given $u^{m-1}\in\Pd_0(\T_h)$, we set the problem: find $u^m\in\Pd_0(\T_h)$ such that
\begin{equation}
  \label{esquema_DG_neuroblasts_unsteady}
  \escalarL{\delta_t u^m}{\bu}+a_\Lambda(u^m,\bu)=L_\Lambda(\bu), 
  \quad \forall \bu\in\Pd_0(\T_h),
\end{equation}
where, for the parameter set $\Lambda$ defined above, $a_\Lambda(\cdot,\cdot)$ is the bilinear form defined in~(\ref{eq:a_Lambda}) and $L_\Lambda(\cdot)$ is the linear form 
\begin{equation}
\label{eq:L_Lambda2}
L_\Lambda(\bu)\coloneqq \gamma\cdot\escalarL{\mathbbm{1}_{\NZdomain} u^{m-1}}{\bu}.
\end{equation}

\section{Numerical Tests and Adjustment of the Model to Real Data}
\label{sec:experimental-param}

The numerical schemes described above have been used to develop computer programs that solve the PDE model  simulating the migration of neuroblasts in the brain, yielding concrete data and comparing it to experimental data obtained from real rodent brains. The final target is calibrating the model parameters so that its output is matched as closely as possible to the real data.  

\subsection{Computer Implementation in a Realistic Domain} 
\label{sec:implementation.realistic}

\begin{figure}
  \begin{center}
    \begin{tabular}{cc}
      \includegraphics[height=0.14\textheight]{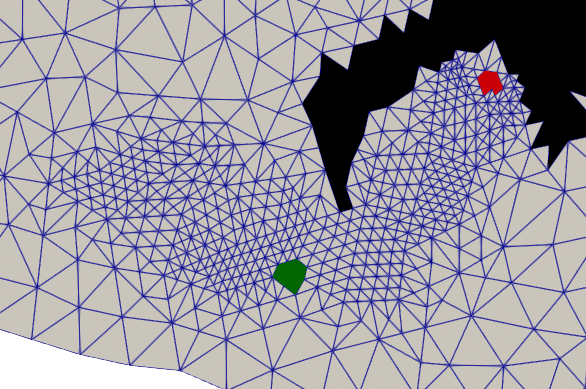}
  &
  \includegraphics[height=0.14\textheight]{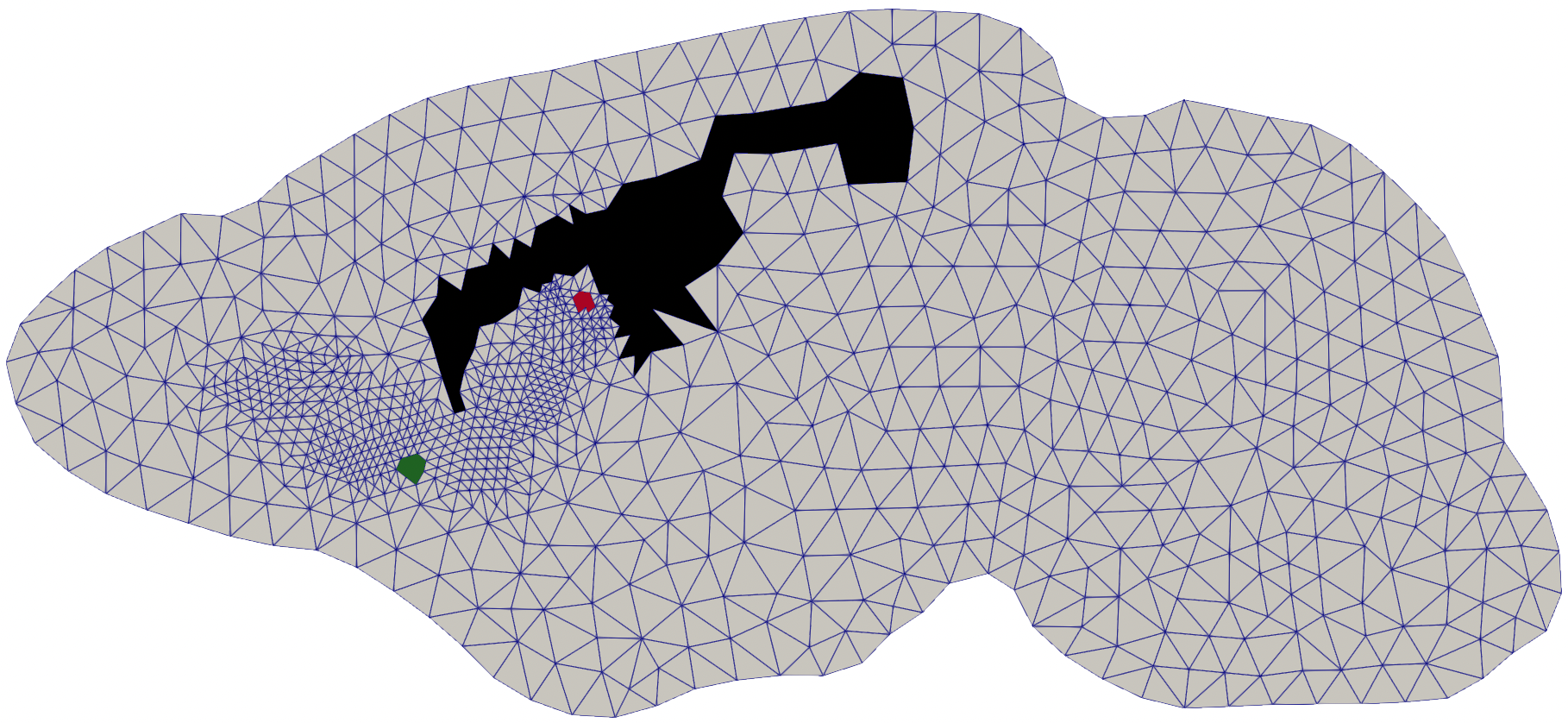}
    \end{tabular}
  \end{center}
    \caption{Right: mesh of a virtual rodent brain. Left: zoom around the RMS. Triangles defining the CC, SVZ and NZ are shown in  black, red and green color, respectively}
  \label{fig:brain-mesh}
\end{figure}

For the computer implementation of the schemes described in Section~\ref{sec:numer-solution}, we have defined a mesh of a virtual rodent brain, see Figure~\ref{fig:brain-mesh}. We used the library \textit{FEniCS}~\cite{alnaes2015fenics} to load and post-process the mesh, define the Galerkin spaces $\Pc_1(\T_h)$ and $\Pd_0(\T_h)$, coding the numerical schemes detailed in Section~\ref{sec:numerical-schemes} and finally storing the solutions for post-processing.

Once we have a suitable mesh, it is post-processed as follows. First, the triangles that are located in the CC, specifically those ones whose barycenter lays on a CC region previously defined, are identified (they are shown in black color in Figure~\ref{fig:brain-mesh}). In the same manner, the triangles laying in the SVZ and in the NZ are identified (they are shown in color red and green, respectively, in Figure~\ref{fig:brain-mesh}). 

This initial mesh is refined  in the area where the RMS is expected to be located and hence where the neuroblast migration occurs.
We identify this area by a preliminary approximation of the initial condition $u_0$ in the initial mesh, which is computed as detailed in Section~\ref{sec:appr-init-cond} for an initial heuristically selected parameter set $\Lambda_0$. Then, those triangles satisfying $u_0>\varepsilon$ (where $\varepsilon$ is a prescribed small constant) are marked and refined leading to a new refined mesh which allows us to compute a more accurate solution for neuroblasts distribution in the region of interest, without significantly increasing the computational effort.

An example of a mesh which was refined using this process can be seen in Figure~\ref{fig:brain-mesh}.
It is the final mesh, $\T_h$, we are using in our numerical tests. It is made of 2003 triangles with size $h \in [3.313\cdot 10^{-3}, 4.814\cdot 10^{-2}]$. 

For the time discretization, we define an uniform partition $t_0<\ldots<t_m<\ldots <t_M$ of the time interval $[0,T]$ with $T=4$ days and constant size $k=t_m-t_{m-1} = 0.04$.

Finally, we code the time advancing schemes detailed in Section~\ref{sec:numerical-schemes}: 
\begin{enumerate}
  \item \emph{Computing $\O$}. For an heuristically computed olfactory bulb shape parameter $\SO$, we compute the $\Pc_1(\T_h)$ solution $\O$ to system~\eqref{esquema_DG_OlfBulb}, where the Gaussian RHS function $f_\O$ depends on $\SO$ and where the anisotropic viscosity $\mu_{O,h}$ introduced in Section~\ref{sec:computing.olf} is defined on those triangles which were marked as located in the CC. For this viscosity constant, we are currently taking $P_C=10^{-30}$, hence in practice the neuroblasts cannot migrate across the corpus callosum. We use standard finite elements to approximate this olfactory bulb function $\O$, which is one of the keys in our model. See Figure~\ref{fig:steady-opt-ob} for an example, where migration of neuroblasts can be anticipated by the isolines of $\O$ (which are orthogonal to $\nabla\O$ and then to the neuroblast path).
  \item \emph{Computing $u^0$}. Given the attraction function $\O$ , we compute a $L^2$ projection of $\nabla\O$ on $\left[\Pc_1(\T_h)\right]^2$. Then, given heuristic positive parameters $\alpha$, $\beta$, $\gamma$, $\X$ (one of which can be eliminated as in Section~\ref{sec:numer-simul-steady}), an initial condition $u^0\in \Pd_0(\T_h)$ can be computed solving the DG method specified in scheme~\eqref{esquema_DG_neuroblasts_steady}. 
    \item \emph{Computing $u^m$, $m> 0$}. Finally, for a new parameter set $\Lambda$ and starting from an initial value $u_0$, for instance the value $u^0\in \Pd_0(\T_h)$ computed above, we approximate at each time step $t_m$, $m=1,\dots,M$, the neuroblast distribution $u^m\in \Pd_0(\T_h)$ by solving the DG scheme~\eqref{esquema_DG_neuroblasts_unsteady}.
\end{enumerate}

\begin{figure}
    \centering
    \includegraphics[width=.8\linewidth]{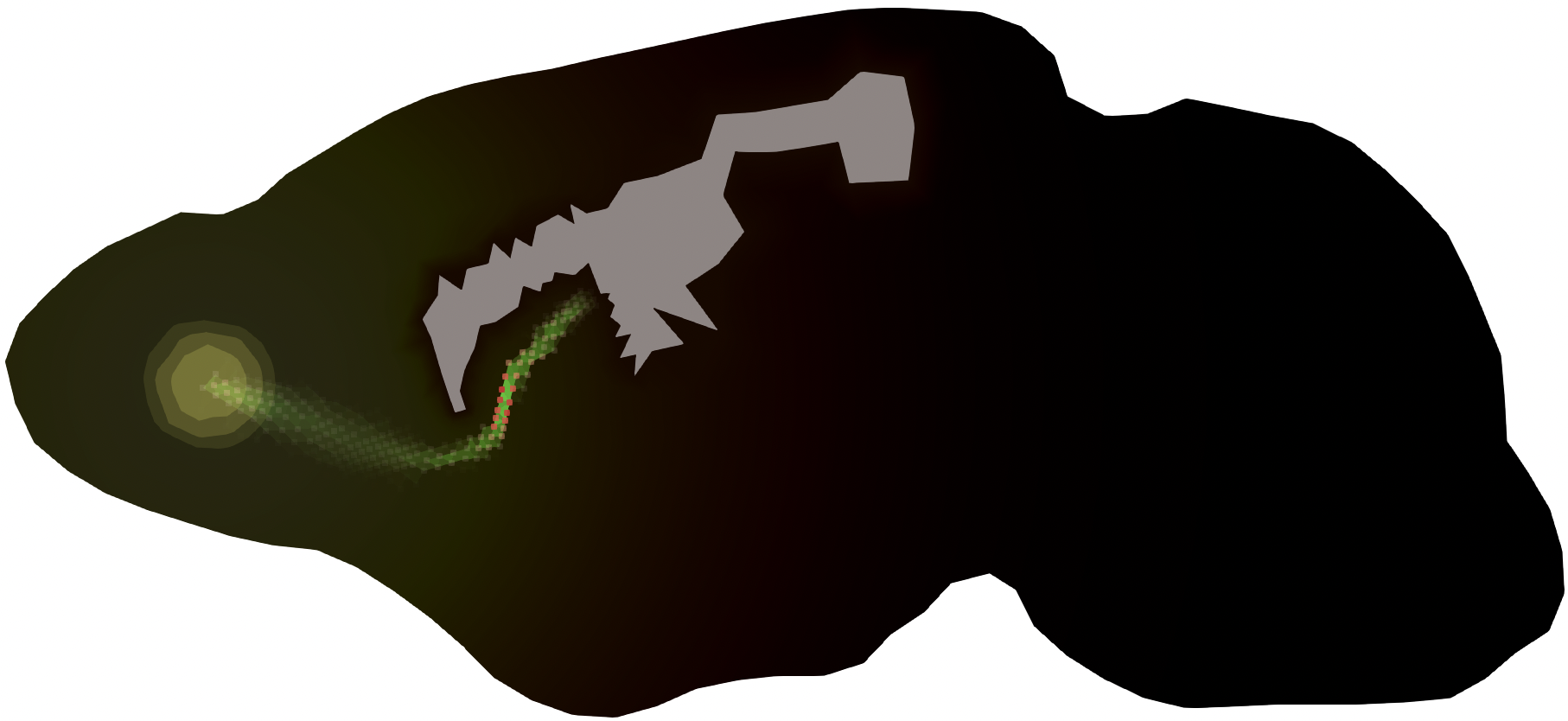}
    \caption{Qualitative behavior of the numerical solution from the stationary model.}
    \label{fig:qualitative}
\end{figure}

As shown in the following sections, the numerical simulations have been really successful from a qualitative point of view, for the parameter set described in Sections \ref{sec:numer-simul-steady} and \ref{sec:numer-simul-evol}. On the one hand, the scheme keeps the positivity of the solution and avoids spurious oscillations. Moreover, as it can be observed in Figure \ref{fig:qualitative}, the solution surrounds the corpus callosum (grey) while it moves along the rostral migratory stream, defined by the chemoattraction, towards the olfactory bulb (yellow).

\subsection{Computing Optimal Parameters}
\label{sec:opt_param}

Although the process described so far is valid and has been successfully tested in qualitative experiments, in what follows we will explore whether it is feasible to adjust the parameters of our model to match the reality. 
Specifically, we describe the process we followed to obtain tuple of parameters $\Lambda$ as in~\eqref{eq:parameters} so that the output of our model is in accordance to the experimental data described next.

\subsubsection*{Errors with respect to Biological Dataset}
In order to analyze whether our model was in agreement with real data obtained from the manual recount of neuroblast migration from the SVZ to the OB of mice during a time course of 7 days, we used the following experimental approach.

Two-month-old male CD1 mice were used for all the experimental procedures. Animals were provided by the Central Animal Experimentation and Production Services (S.E.P.A.) of the Universidad de Cádiz and housed under controlled conditions of temperature (21-23°C) and light (LD 12:12). Mice were fed ad libitum with a commercial diet (AO4 standard maintenance diet, SAFE, Épinay-sur-Orge, France). Experiments were performed strictly following the recommendations given by the Guide for the Care and Use of Laboratory Animals of the European Community Directive 2010/63/UE, and the Spanish regulations (65/2012 and RD53/2013) for the use of laboratory animals. Additionally, all the procedures have been authorized by the Ethics Committee of the “Consejería de Agricultura, Pesca y Desarrollo sostenible” of the Junta de Andalucía (Spain) with the approval numbers 30/03/2016/038 and 04/03/2020/033. 
Mice were placed into labeled cages and divided into five experimental groups, depending on the day of sacrifice (D0, D2, D4, D6, D8; each group n=6). Specifically, the first group (D0) was sacrificed two hours after the last BrdU injection, the second group (D2) was sacrificed at 48 hours, the third (D4) at 96 hours, the fourth (D6) was sacrificed at 144 hours, and finally the D8 group was sacrificed at 192 hours after the start of the experiment. All mice received three intraperitoneal injections of bromodeoxyuridine (BrdU, 120mg/kg) separated by two-hour intervals. Finally, mice were sacrificed and tissue was processed. 
For sacrifice, animals were deeply anesthetized with a lethal dose of 50 mg of Dolethal (Vetoquinol, Lure, France) and perfused with 4 per cent paraformaldehyde (PFA) via the ascending aorta. Brains were removed, postfixed overnight with PFA and cryoprotected by immersion in 30 per cent sucrose solution for 24 hours. Then, serial 30 µm-thick sagittal slices were obtained using a cryotome. Immunohistochemistry was performed as previously described (Rabaneda et al., 2008). Primary antibodies used were rat monoclonal anti-BrdU (1:500) from Abcam (Cambridge, UK) and rabbit polyclonal anti-DCX (1:500, Abcam, Cambridge, UK). Secondary antibodies used were donkey anti-rat IgG labeled with AlexaFluor® 594 (1:1000), from Invitrogen (Carlsbad, CA, USA) and donkey anti-rabbit IgG labeled with Alexafluor ® 488 (1:1000), from Invitrogen (Carlsbad, CA, USA). 
To obtain the images, a Zeiss LSM 900 Airyscan 2 confocal microscope was used. Numerous images were taken of each slice to obtain the complete RMS and were reconstructed offline. Reconstructions of the images were imported into Fiji-ImageJ software for manual quantification. Both, the area and the length of the RMS were measured. Finally, reconstructed images were quantified in order to calculate the velocity of cell migration. Thus, three spots were taken as reference: one at the beginning of the RMS (near to the SVZ), one at the middle of the pathway and another at the end (just before the beginning of the olfactory bulb). A square of known area was drawn at each of the points and within each one the area occupied by the RMS was calculated and the number of BrdU+/DCX+ cells was quantified for each group. Then, the distance from the SVZ to the point where the highest Brdu+/DCX+ cell density was found was divided by the day on which it was reached.
This manual counting of neuroblasts leads to a raw dataset consisting of 28 objects with different attributes containing the information for several slices belonging to 19 rodent brains labelled with BrdU and DCX. Among the included data, we collected the number and density of neuroblasts at different regions, the volume of studied regions, the volume of DCX in studied regions, etc.. In particular, the data is tabulated as follows:
\begin{itemize}
  \item The information of the slices are obtained at 4 different moments in time:  $\widehat t_0=8\mbox{h}$ (initial position), $\widehat t_1=48\mbox{h}$, $\widehat t_2=96\mbox{h}$ or $\widehat t_3=144\mbox{h}$. Specifically, we have the information of a total of 8 slices of 5 mice's brains at $\widehat t_0=8\mbox{h}$, 7 slices of 5 brains at $\widehat t_1=48\mbox{h}$, 6 slices of 4 brains at $\widehat t_2=96\mbox{h}$ and 7 slices of 5 brains at $\widehat t_3=144\mbox{h}$.
  \item Each object contains attributes for each of three different rectangle-shaped control regions, $R_1$, $R_2$, $R_3$, which are respectively located in the vicinity of the SVZ, in the central region of the RMS and near the OB.
\end{itemize}


In order to make use of this information, a preliminary data analysis process was carried out in which several objects were categorized as outliers and removed from the dataset. These outliers were classified using the following criteria:
\begin{enumerate}
    \item First, we considered the raw data obtained for each of the slices of the rodents' brains. The density of the neuroblasts in regions where the volume of the DCX occupied less of 40\% of the total volume was discarded.
    \item Then, the average of the data per mouse using the different slices was computed.
    \begin{enumerate}
        \item If the average position of any of the regions $R_1$, $R_2$ or $R_3$ was not accurately adjusted, the mouse was discarded. This position was considered unusual if the distance between the regions lied further than 1.5 times the interquartile range from the interval determined by the first and third quartiles.
        \item The mice whose average density of neuroblasts in any of the different regions at any of the time instants lied further than 1.5 times the interquartile range from the interval determined by the first and third quartiles were dropped.
    \end{enumerate}
\end{enumerate}


As a result of this data-cleaning process, we ended up with the average information (considering the different slices of the brain) of 4 mice per time step -- a total of 12 mouse-averaged data -- that will be used for our study.

Finally, the average of the data for the 4 different mice was computed at different time steps, obtaining the information shown in Tables~\ref{table:real.data.t0} and~\ref{table:data.evolution}. In addition, the standard deviation of the data is shown in Tables~\ref{table:real.data.t0.std} and~\ref{table:data.evolution.std}. Notice that the data for the region $R_1$ present several difficulties. For instance, at $\widehat t_0=8h$ we only have one object and hence the standard deviation is 0 in this case in Table~\ref{table:real.data.t0.std}. On the other hand, the standard deviation at the later instants $\widehat t_1 = 48h$ and $\widehat t_2 = 96h$ is huge.
\begin{table}[ht]
  \centering
  \begin{tabular}{lccc}
     ($\times 10^4$)& $R_1$ (SVZ) & $R_2$ (RMS) & $R_3$ (OB) \\ \noalign{\smallskip}
    \hline\noalign{\medskip}
    $\widehat t_0 = 8 h$ & $13.671$ & $8.3749$ & $2.9106$
  \end{tabular}
  \caption{Real averaged neuroblast density (neuroblasts/mm$^3$) in three brain regions at the initial time step}
  \label{table:real.data.t0}
\end{table}
\begin{table}[ht]
  \centering
  \begin{tabular}{lccc}
     ($\times 10^4$)& $R_1$ (SVZ) & $R_2$ (RMS) & $R_3$ (OB) \\ \noalign{\smallskip}
    \hline \noalign{\medskip}
    $\widehat t_1 = 48h$ & $7.3009$ & $7.5293$ & $9.9146$\\
    \noalign{\smallskip}
    $\widehat t_2 = 96h$ & $6.7128$ & $5.6462$ & $8.9673$ \\
    \noalign{\smallskip}
    $\widehat t_3 = 144h$ & $1.1033 $ & $3.7577$ & $4.0688 $ 
  \end{tabular}
  \caption{Real averaged neuroblast density (neuroblasts/mm$^3$) in three brain regions at the three later time steps}
  \label{table:data.evolution}
\end{table}

\begin{table}[ht]
  \centering
  \begin{tabular}{lccc}
     ($\times 10^4$)& $R_1$ (SVZ) & $R_2$ (RMS) & $R_3$ (OB) \\ \noalign{\smallskip}
    \hline\noalign{\medskip}
    $\widehat t_0 = 8 h$ & 0 & $3.9534$ & $0.6785$
  \end{tabular}
  \caption{Standard deviation of neuroblast density (neuroblasts/mm$^3$) in three brain regions at the initial time step}
  \label{table:real.data.t0.std}
\end{table}
\begin{table}[ht]
  \centering
  \begin{tabular}{lccc}
     ($\times 10^4$)& $R_1$ (SVZ) & $R_2$ (RMS) & $R_3$ (OB) \\ \noalign{\smallskip}
    \hline \noalign{\medskip}
    $\widehat t_1 = 48h$ & $8.2975$ & $2.7164$ & $1.8624$\\
    \noalign{\smallskip}
    $\widehat t_2 = 96h$ & $11.191$ & $3.3942$ & $6.4012$ \\
    \noalign{\smallskip}
    $\widehat t_3 = 144h$ & $0.8128$ & $0.9098$ & $2.1148$ 
  \end{tabular}
  \caption{Standard deviation of neuroblast density (neuroblasts/mm$^3$) in three brain regions at the three later time steps}
  \label{table:data.evolution.std}
\end{table}


Our purpose is to check whether our model is capable to reproduce, quantitative and qualitatively, this behavior of the neuroblasts migration process.
For that, let  $\widehat u_i^m$ denote the averaged real density of neuroblasts in the brain region $R_i$, at the time step $\widehat t_m$ (Tables~\ref{table:real.data.t0} and~\ref{table:data.evolution}).
We pretend to calculate an initial condition $u^0$ under a certain set of parameters $
\Lambda_0=(\X_0, \SO_0, \alpha_0, \beta_0, \gamma_0)\in\mathbb{R}^5_+,
$ and a second parameter set $\Lambda = (\alpha, \beta, \gamma, \X, \SO)$ such that the related solution to our evolution model minimizes the difference with respect to $\widehat u_i^m$.

More in detail: let us assume  (maybe adjusting the time step) that for each time $\widehat t_k$, $k\in\{0,1,2,3\}$, one can pick some time step $m=m_k\in\mathbb N$ such that $t_{m_k}=\widehat t_k$. 
Let us consider some balls $B_i$, which will be heuristically located in the computational domain $\Omega$ with the aim of reproducing the results in the actual regions $R_i$, $i\in\{1,2,3\}$.  
At the time  $\widehat t_k$, let $u^k_i$ be an approximation of the mean $\int_{B_i}u^{m_k}\,dx/|B_i|$, where $u^{m_k} = u^{m_k}(\Lambda)$ is the approximate solution to equation~\eqref{esquema_DG_neuroblasts_unsteady} for the parameter set $\Lambda$.
We consider the following error expression at time $\widehat t_k\in [0,T]$:
\begin{equation}
\label{error.time.tm}
E^k=E^k(\Lambda)=\frac{1}{N_x}\sum_{i=1}^{N_x}
\frac{1}{\left(\widehat u^{k}_{i}\right)^2}
(u^{k}_{i}  - \widehat{u}^k_{i})^2, 
\end{equation}
where $N_x=3$ is the number of brain regions/computational balls and $k\in\{0,1,2,3\}$.
Notice that $E^0$, the error for $k=0$, shall be normally computed from a different parameter set $\Lambda_0$, as it will be related to the initial condition $u^0$ obtained from the steady equation~\eqref{esquema_DG_neuroblasts_steady}. Hence, it is not used for the global error expression, which is defined as follows, for $N_T=3$:
\begin{equation}
\label{error}
E=E(\Lambda)=
\frac{1}{N_T} \sum_{k=1}^{N_T} E^k =
\frac{1}{N_T\, N_x}\sum_{k=1}^{N_T} \sum_{i=1}^{N_x}
\frac{1}{\left(\widehat u^k_{i}\right)^2}
(u^k_{i}  - \widehat{u}^k_{i})^2.
\end{equation}

With all of the above, we are going to design a process which is split in two steps: 
\begin{enumerate}
  \item Compute a  parameter set $\Lambda_0$ for minimizing the initial error $E^0(\Lambda_0)$, obtaining an initial state, $u_0$, 
\item  Calculate a second parameter set $\Lambda$ for the evolution neuroblast model, obtaining the final solution $u^m$ by minimization of $E(\Lambda)$. 
\end{enumerate}

\subsubsection*{Optimization Process}

For calibrating the parameters, we mainly use the \textsc{L-BFGS-B} (limited\--memory Broyden\-–Fletcher\-–Goldfarb\-–Shanno) method, a popular algorithm for parameter identification in machine learning which falls in the family of quasi-Newton iterative methods~\cite{fletcher2000practical}. 
In particular, we identify the parameters minimizing both the initial error $E^0(\Lambda_0)$ and the evolution error $E(\Lambda)$ with the aid of the Python \texttt{optimize} library~\cite{SciPy-NMeth}, which contains an implementation of the \textsc{L-BFGS-B} method which in addition, allows us to set the bound constraints for the parameters if needed.

The \textsc{L-BFGS-B} algorithm starts with an initial estimate of the parameters, which we provide with the aid of a non-linear regression technique. In this case, we use the \texttt{sklearn} \texttt{Python} library~\cite{scikit-learn}. To obtain a better approximation for the initial vector of parameters, we apply non-linear regression, based on decision trees (with \texttt{RandomForestRegressor} from \texttt{sklearn.ensemble}). 

Before that, we compute a grid of parameter which are needed to train the regressor, by approximating the solution of the discrete schemes, as detailed in Section~\ref{sec:implementation.realistic}, for a extensive grid of values. 
Specifically, they consists, for each parameter set $\Lambda$ on the values of the errors and their correspondent parameters $E^k_i(\Lambda)$, where the biological data in region $R_i$ is compared to our numerical approximations in the computational balls $B_i$. This data is written in \textit{CVS} files with the  format indicated in Table~\ref{tab:cvs.errors},
\begin{table}
  \begin{center}
    \centering
    \begin{tabular}{cccc}
      \hline\noalign{\smallskip}
      $\widehat t_0$ & $ E_1^0 $ & $ E_2^0 $ & $ E_3^0 $\\\noalign{\smallskip}
      $\widehat t_1$ & $ E_1^1 $ & $ E_2^1 $ & $ E_3^1 $\\\noalign{\smallskip}
      $\widehat t_2$ & $ E_1^2 $ & $ E_2^2 $ & $ E_3^2 $\\\noalign{\smallskip}
      $\widehat t_3$ & $ E_1^3 $ & $ E_2^3 $ & $ E_3^3 $\\\noalign{\smallskip}\hline
    \end{tabular}
  \end{center}
  \caption{Format of errors stored in CVS files with the aim of further processing for parameter fitting}\label{tab:cvs.errors}
\end{table}
where $E^m_i=E^m_i(\Lambda)$ is given in~\eqref{error.time.tm}.
The regressor, fed by this data can thus provide a 
first approximation of the parameters. 
More specifically, the steady problem, used to obtain the initial condition, only requires the first row of the table, which corresponds to the initial time. Then, the whole table is used to find the approximated parameters of the time evolution problem. After heuristic revision, they will be the starting point to our final 
\textsc{L-BFGS-B} approximation.

%
%
\subsection{Resulting Parameters for Initial Condition}
\label{sec:numer-simul-steady}

We are now in a position to calculate the appropriate parameters to fit our scheme to the experimental data.
The first step is to find the parameters for the steady problem fitting the initial condition with data at the initial time. In this case, to reduce the number of parameter combinations and then the computational cost, we are going to eliminate parameter $\X$ by dividing the equation by it. This parameter is a good choice for division, because we assume (and our numerical tests confirm this assumption) that it cannot be too small. Then our steady equation is written as:
%
\begin{equation}
  \label{eq:rescaled-stationary}
\aupw{\mmu}{u}{\bv} + \frac{\alpha_0}{\X_0} \escalarL{u}{\bv} - \frac{\beta_0}{\X_0} \escalarL{\mathbbm{1}_{\SVZdomain}}{\bv} =  \frac{\gamma_0}{\X_0} \escalarL{\mathbbm{1}_{\NZdomain}}{\bv}. 
\end{equation}
The bilinear form $\aupw{\mmu}{\cdot}{\dot}$ is defined in~\eqref{eq:aupwFinal} and $\mmu$ is a $P^1$--continuous projection of $\nabla\O$, where $\O$ is an approximation of the solution of equation~\eqref{esquema_DG_OlfBulb} depending on the shape parameter $\SO$ trough the source term for $f_{\O}$, according to~\eqref{source}.

Let us denote the rescaled parameters as 
$$ \alpha_0'=\frac{\alpha_0}{\X_0} , \quad \beta_0'=\frac{\beta_0}{\X_0},\quad \gamma_0'=\frac{\gamma_0}{\X_0}. 
$$ 
Then we focus on finding useful grid of values for these parameters, where we get under and over approximation errors between the real values and the numerical ones in the considered regions.
Taking into account the positivity of the parameters and some previous qualitative numerical tests we select the following intervals: 
$$
\alpha_0' \in [1, 10^2], 
\quad \beta_0' \in [10, 10^8],
\quad  \gamma_0' \in [10, 10^8].$$
Regarding the parameter $\SO$, it has a different nature to $\alpha_0$, $\beta_0$ and $\gamma_0$, being more directly related to the spread along the domain of the OB. Since we are considering normalized dimensions for the domain, we consider $\SO \in [10^{-4},10] $.
Once these ranges of variation of the parameters have been fixed, we finally build a four dimensional grid, consisting of ten values for each parameter ($10^4$ in total) which are logarithmically scaled.  

Next, for each parameter combination in the grid, we proceed solving equation~\eqref{eq:rescaled-stationary} and compute the error respect the experimental values shown in Table~\ref{table:real.data.t0}, extracted from an average of the  manual recount of neuroblasts in real images.
These errors were stored in CVS files as shown in Table~\ref{tab:cvs.errors}. This way, we have the required data to train the random forest regressor.  

Once finished the training process, the following first approximation of the searched parameters was obtained: $$
\alpha_0' = 7.3575, \qquad \beta_0' = 2.0108\cdot10^3, \qquad \SO = 6.2338 \cdot 10^{-1} \quad \text{ and } \quad \gamma_0'= 1.3847 \cdot10^3. $$
\begin{figure}
    \centering
    \includegraphics[width=.8\linewidth]{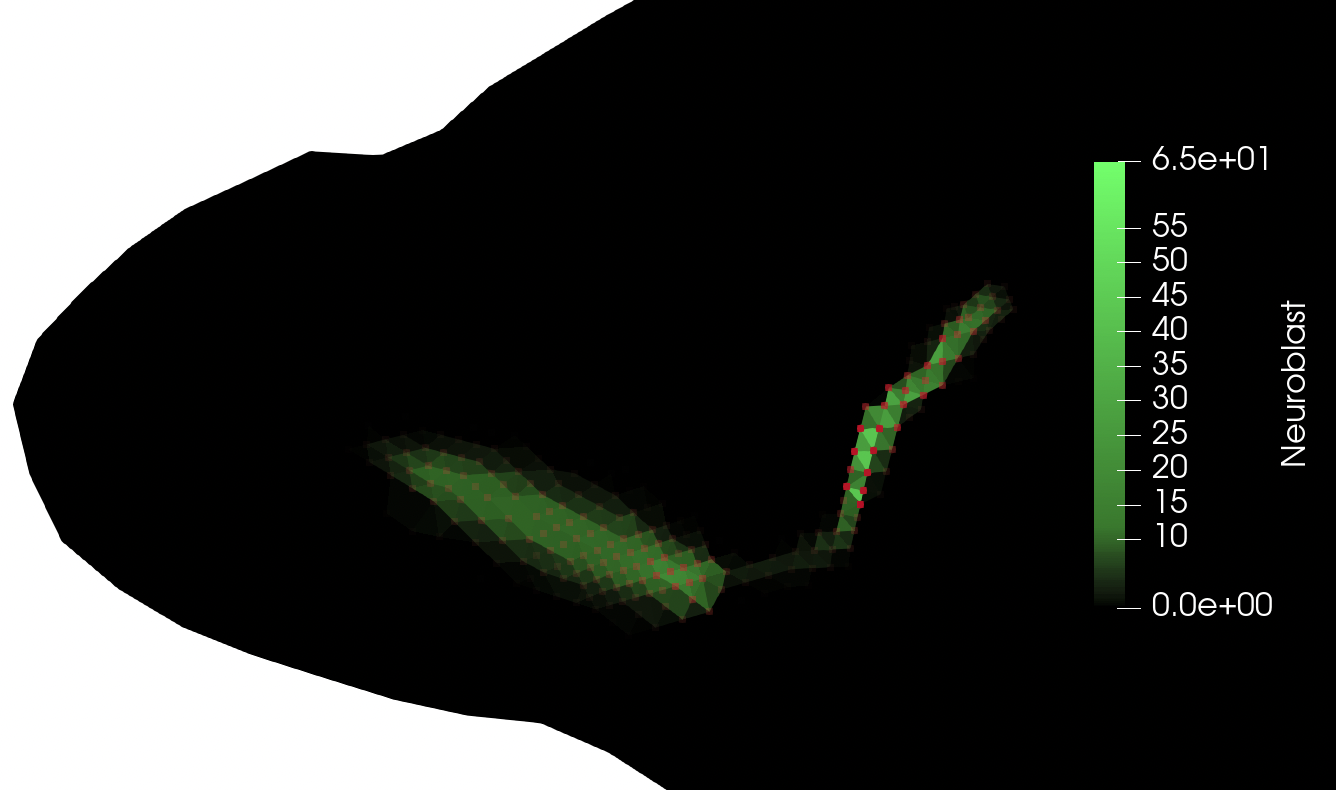}
    \caption{Numerical solution for the initial condition from regression parameters.}
    \label{fig:steady-reg}
\end{figure}
%

\begin{figure}
  \centering
  \includegraphics[width=.8\linewidth]{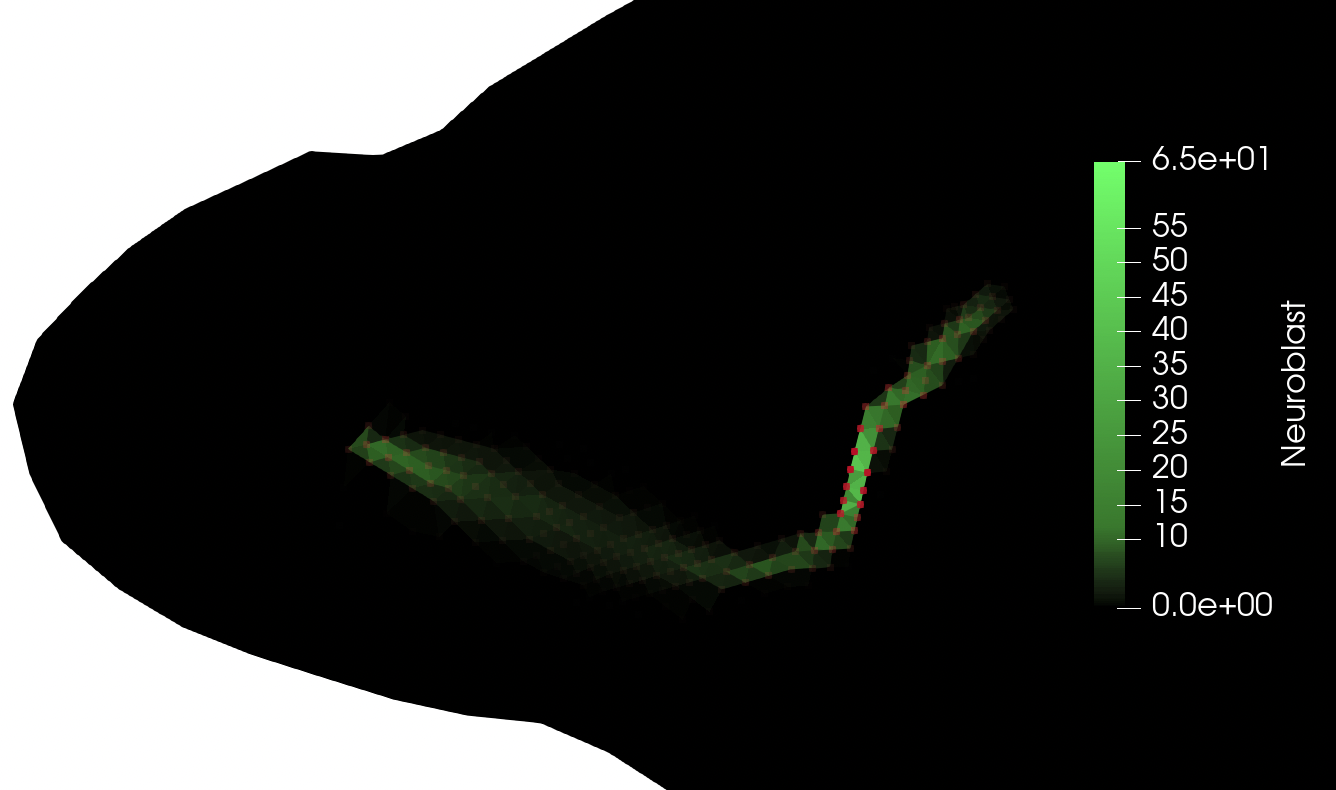}
  \caption{Steady neuroblast density from optimization}
  \label{fig:steady-opt}
\end{figure}

\begin{figure}
  \centering
  \includegraphics[width=.8\linewidth]{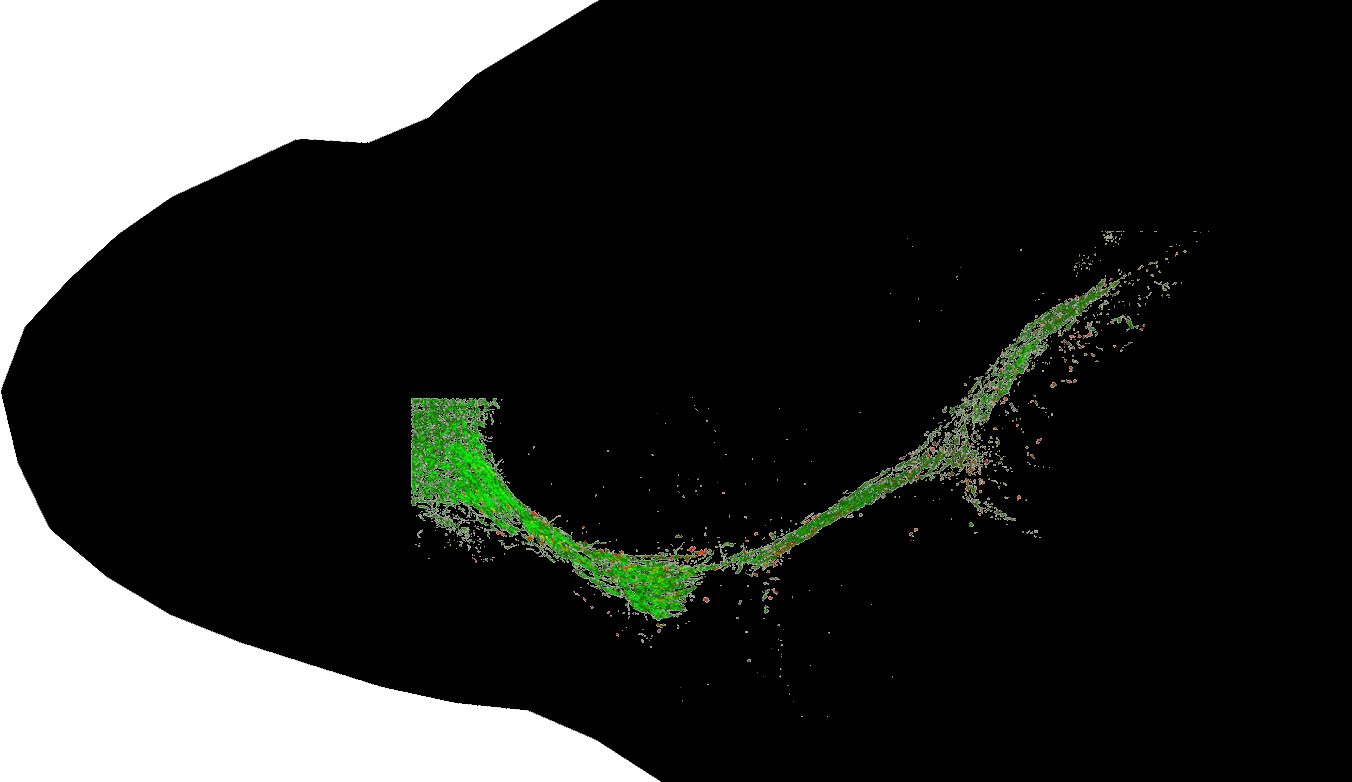}
  \caption{Real steady neuroblast density}
  \label{fig:steady-real}
\end{figure}

These parameters were used as initial values to the \textsc{L-BFGS-B} algorithm, where a maximum number of $10^3$ iterations was imposed in the Python \texttt{optimize} library. Thus, the following parameters were finally obtained:
\begin{equation}
  \label{eq:u0.parameters}
  \alpha_0' = 4.2719, \qquad \beta_0' = 1.7513\cdot10^3, \qquad \SO = 4.0\cdot 10^{-1} \quad \text{ and } \quad \gamma_0'= 0.0. 
\end{equation}
The relative quadratic error obtained was $E = 0.0005$ which means that we have achieved $99.95\%$ of accuracy. The qualitative results, as shown in Figures \ref{fig:steady-opt-ob} and \ref{fig:steady-opt}, agree with the real phenomena from Figure \ref{fig:steady-real}. Even more, the density of neuroblasts computationally obtained in each computational ball for these parameters, which are shown in Table~\ref{table:data.u0}, can be considered as an accurate result compared to the biological data in Table~\ref{table:real.data.t0}, especially if the standard deviation of the experimental data is taken into consideration (see Tables~\ref{table:real.data.t0.std} and~\ref{table:data.evolution.std}).
\begin{table}[ht]
  \centering
  \begin{tabular}{lccc}
     ($\times 10^4$) & $R_1$ (SVZ) & $R_2$ (RMS) & $R_3$ (OB) \\ \noalign{\smallskip}
    \hline \noalign{\medskip}
    $\widehat t_0 = 8 h$ & $13.554$ & $8.6322$ & $2.8413$
  \end{tabular}
  \caption{Data (neuroblasts/mm$^3$) produced by our model for the initial condition}
  \label{table:data.u0}
\end{table}

\subsection{Resulting Parameters for the Neuroblast Evolution Model}
\label{sec:numer-simul-evol}

Now, starting from the initial solution $u^0$ computed from the coefficients obtained in previous section, we are going to find the parameters which make the evolution model fit the real data as much as possible. Unlike the stationary case, we cannot divide by $\X$ to reduce the parameter set due to the additional time derivative term.
 But now we can consider the following fact: we are modeling those neuroblasts that, at the initial time, were marked with bromodeoxyuridine (BrdU), which are represented as red spots in Figure~\ref{fig:neuroblast-real}. Since no new neuroblasts will be marked in later times, we can set $\beta=0$, neglecting the neuroblast birth term in the SVZ. 

So we have again four parameters to fit in our evolution numerical scheme, which can be expanded as:
\begin{equation}
\escalarL{\delta_tu^m}{\bu} + \X \aupw{\nabla \O}{u^m}{\bu} + \alpha \escalarL{u^m}{\bu} = \gamma \escalarL{\mathbbm{1}_{\NZdomain} \,u^{m-1}}{\bu}.
\end{equation}
Now let us consider the intervals $\alpha, \X \in [10^{-4},1]$ and $\gamma \in [10^{-3},10]$, whereas the parameter $\sigma$ is fixed to the value obtained for the stationary problem in the previous section, $\sigma=4.0\cdot 10^{-1}$, see~\eqref{eq:u0.parameters}. In fact, the latter parameter $\sigma$ can be considered as a key physiological constant modeling the spread of the olfactory bulb, therefore, once obtained a value that determines its form, it is not supposed to change.

Now, we build a three dimensional grid with 5 values for each parameter at each time step and we train the random forest regressor as specified in the previous Section~\ref{sec:numer-simul-steady} so that the solution fits the experimental data values shown in Table~\ref{table:data.evolution}. As a result, we obtain a first approximation of the searched parameters:
$$ \alpha = 3.67\cdot10^{-2}, \qquad \X = 4.87 \cdot 10^{-2} 
\quad \text{ and } \quad \gamma = 9.73. $$

\begin{figure}
  \centering
  \includegraphics[width=.8\linewidth]{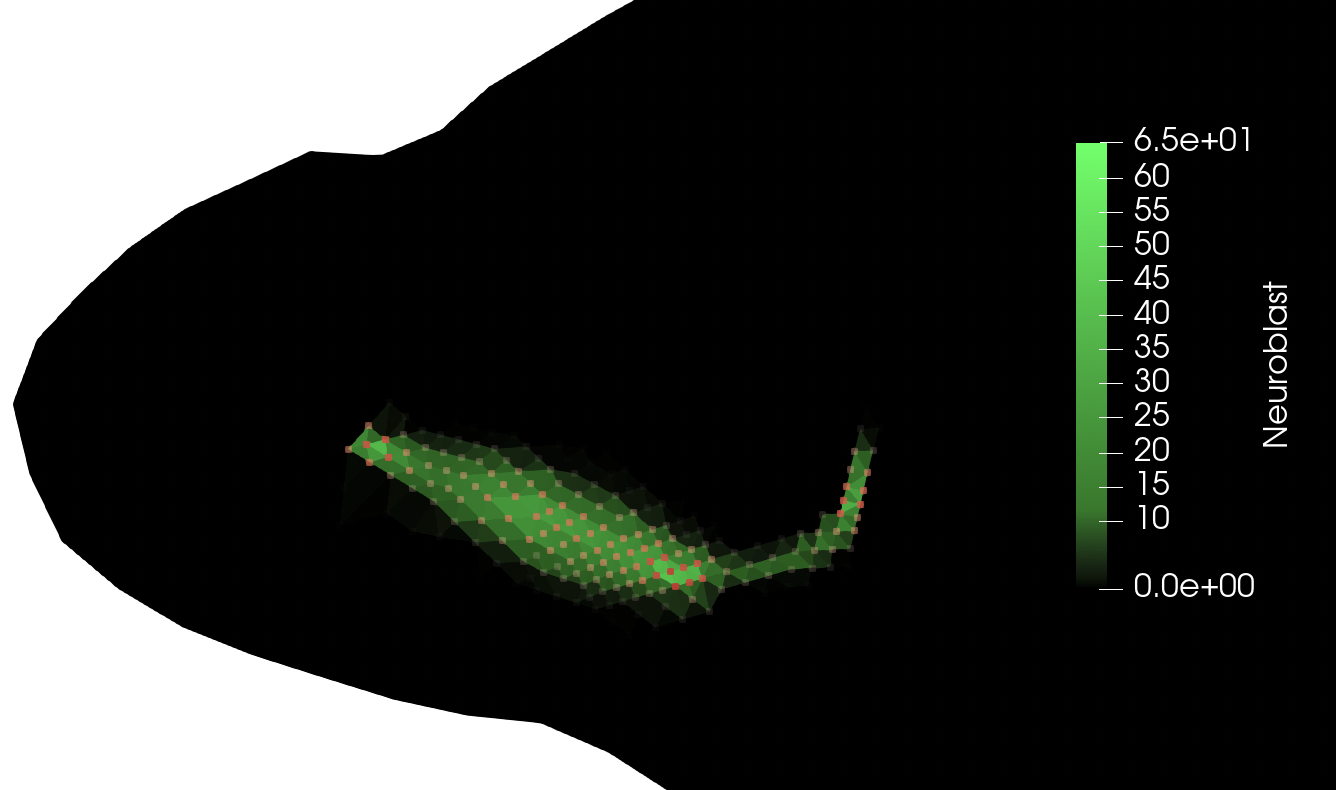}
  \caption{Evolution of neuroblast density after two days (at $\widehat t_1$) for the optimal parameters}
  \label{fig:ev-opt-2}
\end{figure}

\begin{figure}
  \centering
  \includegraphics[width=.8\linewidth]{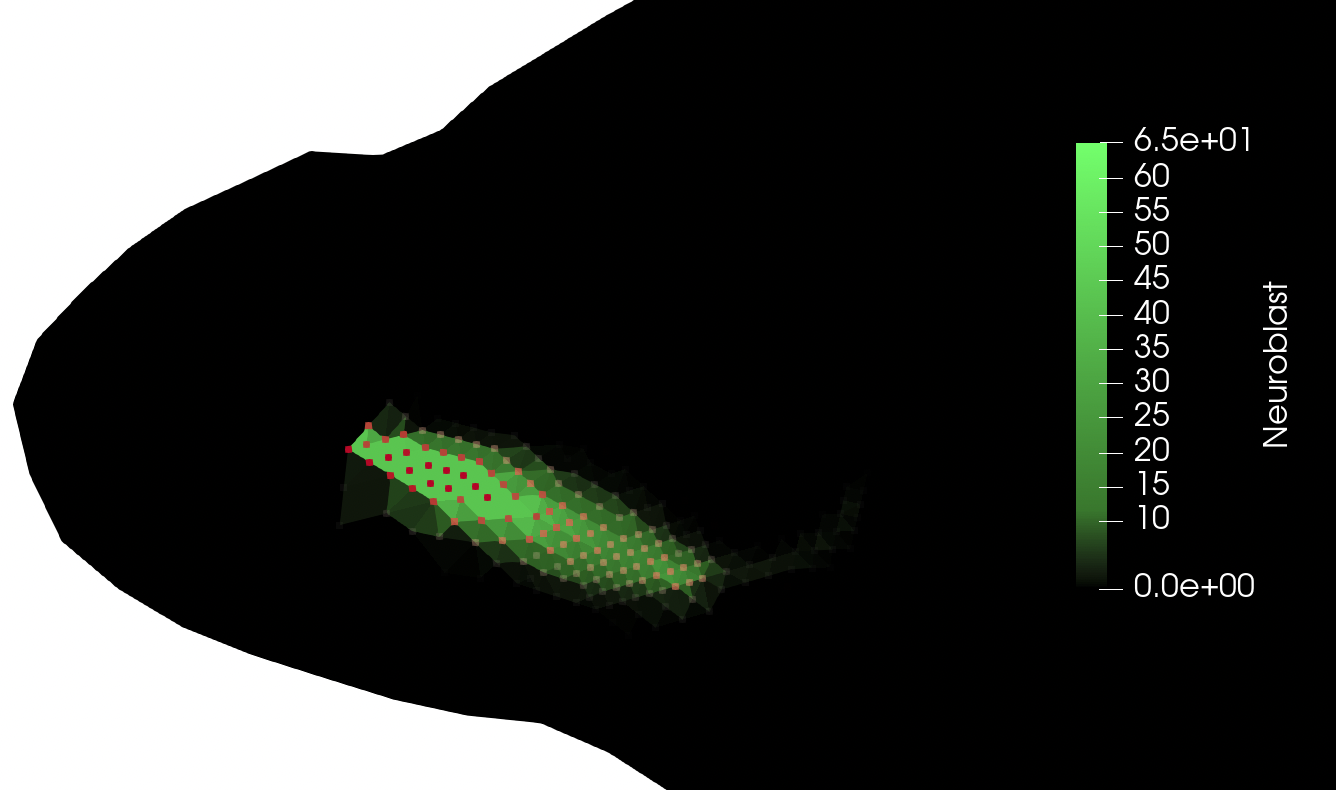}
  \caption{Evolution of neuroblast density after four days (at $\widehat t_2$) for the optimal parameters}
  \label{fig:ev-opt-4}
\end{figure}

\begin{figure}
  \centering
  \includegraphics[width=.8\linewidth]{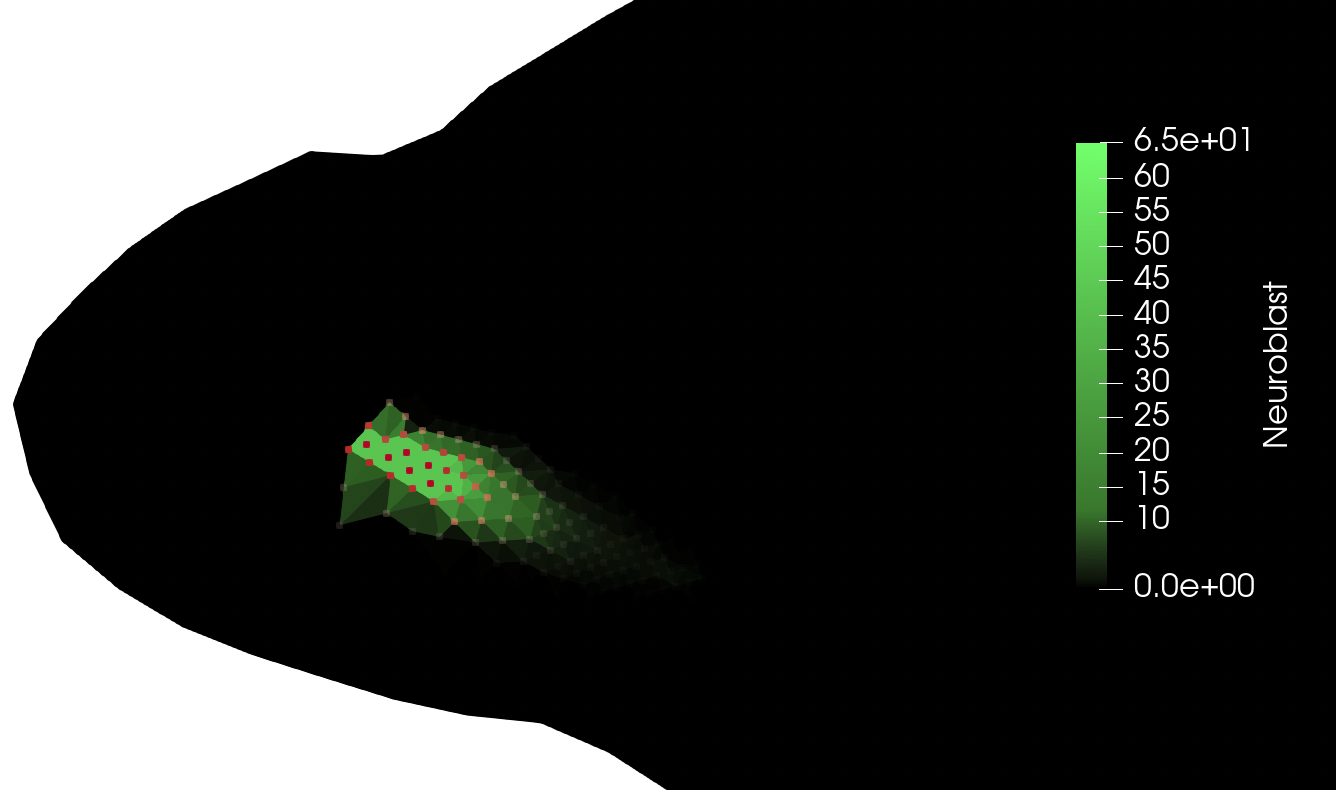}
  \caption{Evolution of neuroblast density after six days (at $\widehat t_3$) for the optimal parameters}
  \label{fig:ev-opt-6}
\end{figure}

Finally, we apply the L-BFGS-B method to optimize the error using the previous values of the parameters as initial guess and we obtain:
$$ \alpha = 1.951\cdot10^{-1}, \qquad \X = 2.241 \cdot 10^{-2} 
\quad \text{ and } \quad \gamma = 3.548. $$
The total relative quadratic error obtained  with respect to the data in Table~\ref{tab:data.evolution.model}, adding the error in each of the time steps $\widehat t_1$, $\widehat t_2$ and $\widehat t_3$, is: $E = 0.38$ ($62\%$ of accuracy). We have represented the neuroblasts migration process for the optimal parameters in Figures~\ref{fig:ev-opt-2}, \ref{fig:ev-opt-4} and~\ref{fig:ev-opt-6}.

If we also consider the approximation of the initial condition, the relative quadratic error made by our simulation of the entire migration process with respect to the real data in Tables~\ref{table:real.data.t0} and~\ref{table:data.evolution} is: $E = 0.28$ ($72\%$ of accuracy). The mean of the density of the neuroblasts in each position $B_i$ and time $\widehat t_k$ is shown in Table~\ref{tab:data.evolution.model}.
\begin{table}[ht]
\centering
  \begin{tabular}{lccc}
     ($\times 10^4$) & $R_1$ (SVZ) & $R_2$ (RMS) & $R_3$ (OB) \\ \noalign{\smallskip}
    \hline \noalign{\medskip}
$\widehat t_0 = 8 h$ & $13.554$ & $8.6322$ & $2.8414$\\\noalign{\smallskip}
$\widehat t_1 =48h$ & $0.0$ & $7.0802$ & $5.5205$\\\noalign{\smallskip}
$\widehat t_2 =96h$ & $0.0$ & $5.7271$ & $6.3969$ \\\noalign{\smallskip}
$\widehat t_3 =144h$ & $0.0$ & $3.1629$ & $5.2144 $
\end{tabular}
\caption{Data (neuroblasts/mm$^3$) produced by our evolution model.}\label{tab:data.evolution.model}
\end{table}
%
%

\section{Conclusions}

In this paper we developed a mathematical PDE model for neuroblast migration in the rodent brain in which migration is mainly determined by a chemoatraction force towards the olfactory bulb (OB). This chemoattraction is determined by the gradient of a very sophisticated function which takes into account the anisotropic morphology of the rodent brain, specifically the presence of a zone, the corpus callosum, where the mobility of neouroblasts is drastically reduced. 
This way, our model can represent the evolution of the density of neuroblasts in realistic two-dimensional brain domains, following a path that corresponds to the rostral migratory stream (RMS) which has been extensively described in the literature.

Even more, we introduced some recent upwind DG numerical schemes for numerically solving the chemotaxis equations. We showed that these schemes are perfectly adapted to our model, both for computing an initial neuroblast distribution as solution of the steady-state equations and, after that, for solving the neuroblast evolution equations.

We have shown some computational tests where we observed that the results obtained from our model are consistent with the real distribution of neuroblasts in rodent brains. And moreover, we carried out some tests where we qualitatively compared the results with respect to a dataset that we previously had obtained from real mouses. In these tests we managed to tune the parameters of our PDE model, so that the results of this model fits the experimental data reasonably well.

Thus, we found that chemoattraction toward specific OB functions is the right way to model the migration of neuroblast in realistic domains. Despite that, of course our model still have room for improvement. It should be studied whether the addition of new terms to the PDEs can improve the fit of the model to the data. And also we should test our model with other datasets, of even higher quality if possible.

In any case, we believe that our preliminary model is a good starting point for a final validated neuroblast model, which should be able to describe neuroblast migration in the adult brain. This is a field that is basically unexplored and that would be a starting point for more general models. For instance, it could be extended to investigate migration of neuroblasts toward brain injuries. And of course to three-dimensional models, which would not pose a great difficulty from a conceptual point of view.

\section*{Acknowledgments}
This publication is part of the R+D+i (PID2022-142418OB-C21) grant funded by MICIU/AEI/10.13039/501100011033 and by FEDER/UE. PPG was granted a fellowship from Ministerio de Universidaddes, Spain (FPU20/00183). DAS has been supported by a UCA FPU contract (UCA/REC14 VPCT/2020) funded by Universidad de Cádiz and by a Graduate Scholarship funded by the University of Tennessee at Chattanooga, USA. FGG and NOR where financed by the project PAIDI 2020 ``Nuevos retos en el estudio de procesos biológicos mediante ecuaciones en derivada parciales" (P20\_01120). DAS, NOR and JRRG were supported by the FQM-315 research group of the Junta de Andalucía, Spain. 


\clearpage
\bibliographystyle{alpha}
\bibliography{references.bib}

\end{document}